\documentclass[12pt]{article}
\usepackage{amsmath}
\usepackage{amssymb}

\usepackage[cp1251]{inputenc}
\usepackage[russian]{babel}
\newcommand {\re} {Re \hskip 1.5pt}
\newcommand{\liml}{\lim\limits_{t \to + \infty}}

\newcommand{\Supl}[1]{\sup\limits_{#1}}
\newcommand{\expb}[1]{\exp\biggl\{#1\biggr\}}
\newcommand{\il}[2]{\int\limits_{#1}^{#2}}
\newcommand{\ilp}[1]{\int\limits_{#1}^{+\infty}}

\newcommand{\ph}{\phantom{a}}
\newcommand{\phh}{\phantom{aaa}}

\newcommand{\sist}[2]{\left\{
\begin{array}{l}
{#1}\\
\ph\\
{#2}
\end{array}
\right.}

\begin{document}
MSC 34C99

\vskip 10pt
{\bf \centerline {Properties of solutions of  quaternionic}

\centerline{Riccati equations}

\vskip 20pt \centerline {G. A. Grigorian}}
\centerline{\it Institute  of Mathematics NAS of Armenia}
\centerline{\it E -mail: mathphys2@instmath.sci.am}
\vskip 20 pt

Abstract. In this paper we study properties of regular solutions of quaternionic Riccati equations. The obtained results we use for study of the asymptotic behavior of solutions of two first-order linear quaternionic ordinary differential equations.

\vskip 20pt
Key words: quaternions, the matrix representation of quaternions,  quaternionic Riccati equations, regular, normal and extremal solutions of Riccati equations, normal, irreconci-\linebreak lable, sub extremal and super extremal systems, principal and non principal solutions.

\vskip 20pt
{\bf 1. Introduction}. Let  $a(t),\phantom{a} b(t),\phantom{a} c(t)$ and $d(t)$  be quaternionic-valued  continuous functions on $[t_0,+\infty)$, i.e.:
$a(t)\equiv a_0(t) + i a_1(t) + j a_2(t) + k a_3(t), \phantom{a} b(t)\equiv b_0(t) + i b_1(t) + j b_2(t) + k b_3(t), \phantom{a} c(t)\equiv c_0(t) + i c_1(t) + j c_2(t) + k c_3(t), \phantom{a} d(t)\equiv d_0(t) + i d_1(t) + j d_2(t) + k d_3(t),$
where  $a_n(t), \phantom{a} b_n(t) \phantom{a}, c_n(t), \phantom{a} d_n(t) \phantom{a} (n=\overline{0,3})$ are
 real-valued continuous functions on  $[t_0,+\infty)$, $i,\phantom{a} j, \phantom{a} k$ are the imaginary unities satisfying the conditions
$$
i^2 = j^2 = k^2 = ijk = -1, \phantom{a} ij = - ji = k. \eqno(1.1)
$$
Consider the quaternionic Riccati equation
$$
q' + q a(t) q + b(t) q + q c(t) + d(t) = 0, \phantom{aaa} t\ge t_0. \eqno (1.2)
$$
Particular cases of this equation appear in various problems of mathematics, in particular in problems of mathematical physics (e. g., in the Euler's vorticity dynamics [1], in the Euler's fluid dynamics [2], in the problem of classification of diffeomorphisms of $\mathbb{S}^4$ [3], and in the other ones [4, 5]). A quaternionic-valued function $q = q(t)$, defined on $[t_1,t_2) \ph (t_0\le t_1 < t_2 \le +\infty)$ is called a solution of Eq. (1.2) on $[t_1,t_2)$, if it is continuously differentiable on $[t_1,t_2)$ and satisfies (1.2) on $[t_1,t_2)$. It follows from the general theory of ordinary differential equations that for every $t_1 \ge t_0$ and $\gamma \in \mathbb{H}$ (here and after $\mathbb{H}$ denotes the algebra of quaternions) there exists $t_2 > t_1 \ph (t_2 \le +\infty)$ such that Eq. (1.2) has the unique solution $q(t)$ on $[t_1,t_2)$, satisfying the initial condition $q(t_1) = \gamma$. Thus for every $t_1 \ge t_0$ and $\gamma \in \mathbb{H}$ a solution $q(t)$ of Eq. (1.2) with $q(t_1) = \gamma$ exists or else on some finite interval $[t_1,t_2)$ or else on $[t_1,+\infty)$. In the last case  the solution $q(t)$ we will call a $t_1$-regular  (or simply regular) solution of Eq. (1.2). Notice that some sufficient conditions for existence of regular solutions are obtained in the works [1], [6], [7]. In the real  case  properties  of regular solutions of Eq. (1.2) are studied in [8]   and have  found several  applications (see [9 - 13]). In this paper we study the properties of regular solutions of Eq. (1.2). We use the obtained result to study the asymptotic behavior of solutions of systems of two first-order linear  quaternionic differential equations.

\hskip 10 pt

{\bf 2. Auxiliary propositions}. It is not difficult to verify that there exists a one to one correspondence $q \leftrightarrow Q$ between the quaternions $q = q_0 + i q_1 + j q_2 + k q_3, \ph q_k \in \mathbb{R}, \ph k= \overline{0,3}$ and  the skew symmetric matrices
$$
Q\equiv\left(\begin{array}{rccr}
q_0& q_1& q_2& -q_3 \\
-q_1& q_0 &-q_3 & -q_2\\
-q_2& q_3& q_0& q_1 \\
q_3& q_2& -q_1& q_0
\end{array}
\right),
$$
keeping the arithmetic operations: $q_m \leftrightarrow Q_m, \ph m=1,2 \ph  \Rightarrow q_1 + q_2 \leftrightarrow Q_1 + Q_2, \linebreak q_1 q_2 \leftrightarrow~ Q_1 Q_2, \ph  q_1^{-1} \leftrightarrow Q_1^{-1} \ph (q_1 \ne 0)$.
The matrix $Q$ we will call the symbol of $q$ and will denote by $\widehat{q}$. By $|q|$ we denote the euclidian norm of the vector $q: \ph |q| \equiv \sqrt{q_0^2 + q_1^2 + q_2^2 + q_3^2}$. We also denote $Re\hskip 3pt q \equiv q_0$ - the real part of $q$ and $Im\hskip 3pt q \equiv i q_1 + j q_2 + k q_3$ - the imaginary part of $q$. Finally by $tr \hskip 3pt \widehat{q}$ we denote the trace of $\widehat{q}$.

{\bf Lemma 2.1.} {\it For every quaternion $q$ the equalities
$$
det \hskip 3pt \widehat{q} = |q|^4, \phh tr \hskip 3pt \widehat{q} = 4 Re \hskip 3pt q
$$
are valid.}

Proof. By direct checking.

Let $A(t), \ph B(t), \ph C(t)$ and $D(t)$ be the symbols of $a(t), \ph b(t), \ph c(t)$ and $d(t)$ respectively. Consider the matrix Riccati equation
$$
Y' + Y A(t) Y + B(t) Y + Y C(t) + D(t) = 0, \phh t \ge t_0. \eqno (2.1)
$$
Obviously the solutions $q(t)$ of Eq. (1.2), existing on an interval $[t_1, t_2) \ph (t_0 \le t_1 < t_2 \le +\infty)$ are connected with solutions $Y(t)$ of Eq. (2.1) by relation
$$
\widehat{q(t)} = Y(t), \ph t\in [t_1, t_2). \eqno (2.2)
$$
Let $Y(t)$ be a solution of Eq. (2.1) on $[t_1,t_2)$ Then every solution $Y_1(t)$ of Eq. (2.1) on $[t_1, t_2)$ is connected with $Y(t)$ by the formula (see [14], pp. 139, 140, 158, 159, Theorem~ 6.2)\linebreak
$$
Y_1(t) = Y(t) + [\Phi_Y(t)\Lambda^{-1}(t_1)(I + \Lambda(t_1)\mathcal{M}_Y(t_1,t))\Psi_Y(t)]^{-1},\phh t \in [t_1,t_2),
$$
where $\Phi_Y(t)$ and $\Psi_Y(t)$ are the solutions of the linear matrix equations
$$
\Phi' = [A(t) Y(t) + C(t)]\Phi, \phh t\in [t_1,t_2),
$$
$$
\Psi' = \Psi[B(t) + Y(t) A(t)], \phh t \in [t_1, t_2)
$$
respectively with $\Phi_Y(t_1) = \Psi_Y(t_1) = I, \ph I$ is the identity matrix of dimension $4\times 4$,
$$
\mathcal{M}_Y(t_1,t) \equiv \il{t_1}{t}\Phi_Y^{-1}(\tau)A(\tau)\Psi_Y^{-1}(\tau) d \tau, \phh t\in [t_1,t_2),
$$
$\Lambda(t_1) \equiv Y_1(t_1) - Y(t_1),$ provided $det \Lambda(t_1) \ne 0$. From here we obtain
$$
Y_1(t) = Y(t) +\Psi_Y^{-1}(t) [I + \Lambda(t_1)\mathcal{M}_Y(t_1,t)]^{-1}\Lambda(t_1)\Phi_Y^{-1}(t), \phh t\in [t_1,t_2). \eqno (2.3)
$$
By the    Liouville   formula we have:
$$
det \Phi_Y(t) = \exp\biggl\{\il{t_1}{t} tr [A(\tau) Y(\tau) + C(\tau)] d\tau\biggr\},  \phh t\in [t_1,t_2), \eqno (2.4)
$$
$$
det \Psi_Y(t) = \exp\biggl\{\il{t_1}{t} tr [A(\tau) Y(\tau) + B(\tau)] d\tau\biggr\},  \phh t\in [t_1,t_2), \eqno (2.5)
$$
Let $q(t)$ be a solution of Eq. (1.2) on $[t_1,t_2)$. Then due to (2.2) from (2.3) it follows that for every solotion $q_1(t)$ of Eq. (1.2) on $[t_1,t_2)$ the equality
$$
q_1(t) = q(t) + \psi_q^{-1}(t)[1 + \lambda(t_1)\mu_q(t_1,t)]^{-1} \lambda(t_1)\phi_q^{-1}(t), \phh t\in [t_1,t_2)  \eqno (2.6)
$$
is valid, where $\phi_q(t)$ and $\psi_q(t)$ are the solutions of the linear equations
$$
\phi' = [a(t) q(t) + c(t)]\phi, \phh t\in [t_1,t_2),
$$
$$
\psi' = \psi[b(t) + q(t) a(t)], \phh t\in [t_1,t_2)
$$
respectively with $\phi_q(t_1) = \psi_q(t_1) = 1, \ph \lambda(t_1) \equiv q_1(t_1) - q(t_1),$
$$
\mu_q(t_1,t) \equiv \il{t_1}{t}\phi_q^{-1}(\tau) a(\tau) \psi_q^{-1}(\tau) d \tau, \phh t \in [t_1,t_2).
$$
By (2.3) and Lemma 2.1 from (2.5) and (2.6) we obtain
$$
|\phi_q(t)| = \exp\biggl\{\il{t_1}{t} Re [a(\tau) q(\tau) + c(\tau)] d \tau\biggr\}, \phh t\in [t_1,t_2), \eqno (2.7)
$$
$$
|\psi_q(t)| = \exp\biggl\{\il{t_1}{t} Re [a(\tau) q(\tau) + b(\tau)] d \tau\biggr\}, \phh t\in [t_1,t_2). \eqno (2.8)
$$
Let $q_m(t), \ph m=1,2$ be solutions of Eq. (1.2) on $[t_1,t_2)$. Set: $\lambda_{m,s}(t_1) \equiv q_m(t_1) - q_s(t_1), \linebreak m,s = 1,2.$ By (2.4) we have
$$
a(t)[q_m(t) - q_s(t)] = a(t)\psi_{q_s}^{-1}(t)[1 + \lambda_{m,s}(t_1)\mu_{q_s}(t_1;t)]^{-1}\phi_{q_s}^{-1}(t), \ph t\in [t_1,t_2).
$$
Hence,
$$
[1 + \lambda_{m,s}(t_1)\mu_{q_s}(t_1;t)]' = A_{q_m,q_s}(t_1;t)[1 + \lambda_{m,s}(t_1)\mu_{q_s}(t_1;t)], \phh t\in [t_1,t_2),
$$
where
$$
 A_{q_m,q_s}(t_1;t) \equiv \lambda_{m,s}(t_1)\psi_{q_s}^{-1}(t)[q_m(t) - q_s(t)]\phi_{q_s}^{-1}(t)\lambda_{m,s}^{-1}(t_1), \phh t\in [t_1,t_2), \ph m=1,2.
$$
From here it follows
$$
[I + \widehat{\lambda_{m,s}(t_1)}\widehat{\mu_{q_s}(t_1;t)}]' = \widehat{A_{q_m,q_s}(t_1;t)}[I + \widehat{\lambda_{m,s}(t_1)}\widehat{\mu_{q_s}(t_1;t)}], \phh t\in [t_1,t_2), \ph m=1,2.
$$
By Lemma 2.1 and the Liouville's formula from here we obtain
$$
|1 + \lambda_{m,s}(t_1)\mu_{q_s}(t_1;t)| = \exp\biggl\{\il{t_1}{t} Re [a(\tau) (q_m(\tau) - q_s(\tau))] d\tau\biggr\}, \ph t\in [t_1,t_2), \eqno (2.9)
$$
$m,s = 1,2$. From here we immediately get:
$$
|1 + \lambda_{m,s}(t_1)\mu_{q_s}(t_1;t)| |1 + \lambda_{s,m}(t_1)\mu_{q_m}(t_1;t)|  \equiv 1, \ph t\in [t_1,t_2), \ph m,s =1,2. \eqno (2.10)
$$

{\bf 3. Properties of regular solutions of Eq. (1.2).}

{\bf Definition 3.1.} {\it A $t_1$-regular solution $q(t)$ of Eq. (1.2) is called $t_1$-normal if there exists a neighborhood $U(q(t_1))$ of $q(t_1)$ such that every solution $\widetilde{q}(t)$ of Eq. (1.2) with $\widetilde{q}(t_1) \in U(q(t_1))$ is also $t_1$-regular, otherwise $q(t)$ is called $t_1$-extremal.}

{\bf Definition 3.2.} {\it Eq. (1.2) is called regular if it has at least one regular solution.}

{\bf Remark 3.1.} {\it Since the solutions of Eq. (1.2) are continuously dependent on their initial values every $t_1$-normal ($t_1$-extremal) solution of Eq. (1.2) is also a $t_2$-normal ($t_2$-extremal) solution of Eq. (1.2) for all $t_2 > t_1$. Due to this a $t_1$-normal ($t_1$-extremal) solution of Eq. (1.2) we will just call a normal (a extremal) solution of Eq. (1.2). Note that a $t_2$-normal ($t_2$-extremal) solution of Eq. (1.2) may not be a $t_1$-normal ($t_1$-extremal) solution of Eq.(1.2) if $t_1 < t_2$, because a $t_2$-regular solution of Eq. (1.2) may not be $t_1$-regular for $t_1< t_2$.}

{\bf Theorem 3.1.} {\it If Eq. (1.2) has a $t_1$-regular solution $q(t)$ for some $t_1 \ge t_0$, then it has also another (different from $q(t)$) $t_1$-regular solution.}

Proof. Let $q(t)$ be a $t_1$-regular solution for some $t_1 \ge t_0$. Since $\mu_q(t_1;t)$ is continuously differentiable by $t$ there exists $\gamma \in \mathbb{H} \backslash \{0\}$ such that $\mu_q(t_1;t) \ne \gamma$ for all $t \ge t_0$ ($\mu_q(t_1;t_1) =~ 0$ and the curve $f(t) \equiv \mu_q(t_1;t), \ph t \ge t_1$ is not space filling). Therefore by (2.7) the solution $q_1(t)$ of Eq. (1.2) with $q_1(t_1) = q(t_1) - \frac{1}{\gamma}$ is a $t_1$-regular solution of Eq. (1.2), different from $q(t)$. The theorem is proved.

Denote by $Q(t;t_1;\lambda)$ the general solution of Eq. (1.2) in the region $G_{t_1} \equiv \{(t;q) : \ph t \in~ I_{t_1}(\lambda), \ph  q, \lambda \in\mathbb{H}\},$ where $I_{t_1}$ is the maximum existence interval for the solution $q(t)$ of Eq. (1.2) with $q(t_1) = \lambda$.

{\bf Example 3.1.} {\it Consider the equation
$$
q' + q a(t) q = 0, \phh t \ge -1. \eqno (3.1)
$$
The general solution of this equation in the region $G_0 \cap [-1, +\infty) \times \mathbb{H}$ is given by formula
$$
Q(t;0;\lambda) = \frac{1}{1 + \lambda \il{t_1}{t} a(\tau) d \tau} \lambda, \phh \lambda \in \mathbb{H}, \phh 1 + \lambda \il{t_1}{t} a(\tau) d \tau \ne 0, \phh t \ge t_1. \eqno (3.2)
$$
Assume $a(t)$ has a bounded support. Then from (3.2) is seen that Eq. (3.1) has no $0$-extremal solution, and all its solutions $Q(t,;0;\lambda)$ with  enough small $|\lambda|$ are $0$-normal. If $a(t)$ is a non negative function with an unbounded support and $I_0 \equiv \ilp{0}a(\tau) d \tau < +\infty$ then from (3.2) is seen that the solution $q_0(t) = Q(t;0;-\frac{1}{I_0})$ is $0$-extremal; all the solutions $Q(t;0;\lambda)$ with $\lambda \in \mathbb{H}\backslash (-\infty, - \frac{1}{I_0})$ are $0$-normal and all the solutions $Q(t;0;\lambda)$ with $\lambda \in (-\infty, -\frac{1}{I_0})$ are not $0$-regular. Assume now $\il{0}{t} a(\tau) d \tau = \arctan(\cos t + i \sin t + j \cos \pi t + k \sin \pi t), \ph t \ge 0$. Then from (3.2) is seen that all the solutions $Q(t;0;\lambda)$ with $|\lambda| = \frac{\sqrt{2}}{\pi}$ are $0$-extremal (since the set $\{\frac{1}{\sqrt{2}}(\cos t + i \sin t + j \cos \pi t + k \sin \pi t) : t \ge 0\}$ is everywhere  dense in the unite sphere $\{ q : |q| =1\}$) and all solutions $Q(t;0;\lambda)$ with $|\lambda| < \frac{\sqrt{2}}{\pi}$ are $0$-normal.}

{\bf Example 3.2.} {\it For $u_0 \in \mathbb{H}$ and $0 < r < R < +\infty$ denote $K_{r,R}(u_0) \equiv \{q \in \mathbb{H} : r < |q - u_0| < R\}$ - an annulus in $\mathbb{H}$ with a center $u_0$ and radiuses $r$ and $R$. For any $\varepsilon > 0$ denote $K_{\varepsilon, r,R}(u_0) \equiv \{\xi_1, ..., \xi_m \in K_{r,R}(u_0): \ph if \ph u \in K_{r,R}(u_0) \ph then \ph  there \ph exists \ph s \in~ \{1, ..., m\} \linebreak such \ph that \ph |u - \xi_s| < \varepsilon\}$- a finite $\varepsilon$-net for $K_{r,R}(u_0)$ (here $m$ depends on $\varepsilon$). Consider the sequence of $\frac{1}{2n}$-nets: $\{K_{\frac{1}{2n}, \frac{1}{n}, n}(u_0)\}_{n=1}^{+\infty}$. Let the function $f(t) \equiv \il{0}{t} a(\tau) d \tau, \ph t \ge 0$ has the following properties: $f(t) \ne u_0, \ph t\in [0,1]$; when t varies from $n$ to $n+1 \ph (n=1, 2, ...)$ the curve $f(t)$ crosses all points of $K_{\frac{1}{2n},\frac{1}{n}, n}(u_0)$ (i. e. for every $v \in K_{\frac{1}{2n},\frac{1}{n}, n}(u_0)$ there exists $\zeta_v \in [n,n+1]$ such that $f(\zeta_v) = v); \ph  f(t) \in K_{\frac{1}{2n}, +\infty}(u_0) \ph n =1, 2, .... , \ph t \ge 1.$ From these properties it follows that for every $T\ge 0$ the set $\{f(t) : t\ge T\}$ is everywhere dense in $\mathbb{H}$ and $f(t) \ne u_0, \ph t \ge 0$. Hence from (3.2) it follows that Eq. (3.1) has no $t_1$-normal solutions for all $t_1 \ge 0$ and has at least two extremal solutions: $q_1(t) \equiv 0$ and $q_2(t)$ with $q_2(0) = - \frac{1}{u_0}.$ By analogy using $\frac{1}{2n}$-nets $K_{\frac{1}{2n},\frac{1}{n},n}(u_0;...u_l)\equiv \{\xi_1, ..., \xi_m \in \bigcap\limits_{k=0}^l K_{\frac{1}{n},n}(u_k) : u \in~  \bigcap\limits_{k=0}^l K_{\frac{1}{n},n}(u_k)\Rightarrow \exists s \in \{1, ..., m\} : |u - \xi_s| < \frac{1}{2n}\}$ of the intersections $\bigcap\limits_{k=1}^l K_{\frac{1}{n},n}(u_k)$ in place of $K_{\frac{1}{2n},\frac{1}{n},n}(u_0), \ph n =1, 2, ....$ one can show that there exists a Riccati equation which has no $t_1$-normal solutions and has at least $l+2$ $t_1$-extremal solutions for all $t_1 \ge 0$.}

{\bf Theorem 3.2.} {\it A $t_1$-regular solution $q(t)$ of Eq. (1.2) is $t_1$-normal if and only if $\mu_q(t_1;t)$ is bounded by $t$.}

Proof.  Sufficiency. Set $M\equiv \sup\limits_{t \ge t_1}|\mu_q(t_1;t)|$. Let $q_1(t)$ be a solution of Eq. (1.2) with $|q(t_1) - q_1(t_1)| < \frac{M}{2}.$ Then obviously
$$
1 + (q_1(t_1) - q(t_1)) \mu_q(t_1;t) \ne 0, \phh t \ge t_1.
$$
By (2.7) from here it follows that $q_1(t)$ is $t_1$-normal.

Necessity. Suppose $\mu_q(t_1;t)$ is unbounded by $t$ on $[t_1,+\infty)$. Let then $t_1 < t_2 < ... t_m , ...$ be an infinitely large sequence such that
$$
|\mu_q(t_1;t_n)| \ge n, \phh n =2, 3, ...  \eqno (3.3)
$$
Let $q_n(t), \ph n = 2, 3, ...$ be the solutions of Eq. (1.2) with
$$
q_n(t_1) - q(t_1) = - \mu_q(t_1;t_n) ^{-1}, \phh n =2, 3, ... \eqno (3.4)
$$
Since $q(t)$ is $t_1$-normal there exists $\delta > 0$ such that every solution $\widetilde{q}(t)$ of Eq. (1.2) with $|\widetilde{q}(t_1) - q(t_1)| < \delta$ is $t_1$-regular. Hence from (3.3) and (3.4) it follows that for enough large $n$ the solutions $q_n(t)$ are $t_1$-regular. On the other hand by (2.7) from (3.4) it follows that for enough large $n$ every solution $q_n(t)$ is unbounded in the neighborhood of $t_n$. It means that for enough large $n$ the solutions $q_n(t)$ are not $t_1$-regular. The obtained contradiction completes the proof of the theorem.

By (2.10) from Theorem 3.2 we immediately obtain

{\bf Corollary 3.1.} {\it The following statements are valid:

\noindent
1) any two $t_1$-regular solutions $q_1(t)$ and $q_2(t)$ of Eq. (1.2) are $t_1$-normal if and only if the function
$$
I_{q_1,q_2}(t) \equiv \il{t_1}{t} Re [a(\tau)(q_1(\tau) - q_2(\tau))] d\tau, \phh t \ge t_1
$$
is bounded;

\noindent
2) if $q_N(t)$ and $q_*(t)$ are $t_1$-normal and $t_1$-extremal solutions of Eq. (1.2) respectively then
$$
\limsup\limits_{t \to + \infty} \il{t_1}{t} Re [a(\tau)(q_*(\tau) - q_N(\tau))] d \tau < +\infty,
$$
$$
\liminf\limits_{t \to + \infty} \il{t_1}{t} Re [a(\tau)(q_*(\tau) - q_N(\tau))] d \tau = -\infty;
$$

\noindent
3) if $q_*(t)$ and $q^*(t)$ are $t_1$-extremal solutions of Eq. (1.2) then
$$
\limsup\limits_{t \to + \infty} \il{t_1}{t} Re [a(\tau)(q_*(\tau) - q^*(\tau))] d \tau = +\infty,
$$
$$
\liminf\limits_{t \to + \infty} \il{t_1}{t} Re [a(\tau)(q_*(\tau) - q^*(\tau))] d \tau = -\infty.
$$}
$\phantom{aaaaaaaaaaaaaaaaaaaaaaaaaaaaaaaaaaaaaaaaaaaaaaaaaaaaaaaaaaaaaaaa} \blacksquare$

{\bf Definition 3.3.} {\it A regular Eq. (1.2) is called normal if it has no extremal solutions.}

{\bf Definition 3.4.} {\it A regular Eq. (1.2)  is called  irreconcilable  if its every regular solution  is extremal.}

{\bf Definition 3.5}. {\it A regular  Eq. (1.2) is called sub extremal  if it has only one extremal solution.}

{\bf Definition 3.6.} {\it A regular Eq. (1.2) is called super extremal if it has at least two extremal solutions and normal solutions.}

From Definitions 3.3 - 3.6 is seen that every regular Eq. (1.2) is or else normal or else irreconcilable or else sub extremal or else super extremal. The examples, illustrated above, show that all these types of equations exist.

For any $t_1$-regular solution $q(t)$ of Eq. (1.2) set
$$
\nu_q(t) \equiv \ilp{t}\phi_q^{-1}(\tau) a(\tau) \psi_q^{-1}(\tau) d \tau, \phh t \ge t_1,
$$
where   $\phi_q(t)$ and $\psi_q(t)$ are the solutions of the linear equations
$$
\phi' = [a(t) q(t) + c(t)]\phi, \phh t\ge t_1.
$$
$$
\psi' = \psi[b(t) + q(t) a(t)], \phh t\ge t_1
$$
respectively with $\phi_q(t_1) = \psi_q(t_1) = 1.$

{\bf Theorem 3.3.} {\it Let $q_0(t)$ be a $t_1$-regular solution of Eq. (1.2) such that the integral $\nu_{q_0}(t_1)$ is convergent. Then in order that Eq. (1.2) has a $t_1$-extremal solution it is necessary and sufficient that $\nu_{q_0}(t) \ne 0, \ph t \ge t_1$. If this condition is satisfied then:

\noindent
1) the unique $t_1$-extremal solution $q_*(t)$ of Eq. (1.2) is given by the formula
$$
q_*(t) = q_0(t) - \frac{1}{\nu_{q_0}(t)}, \phh t \ge t_1; \eqno (3.5)
$$

\noindent
2) for all $t_1$-normal solutions $q(t)$ of Eq. (1.2) and  only for them the integrals $\nu_q(t)$ converge for all $t \ge t_1$ and $\nu_q(t) \ne 0, \ph t \ge t_1;$

\noindent
3) for all $t\ge t_1$
$$
\nu_{q_*}(t) = \infty; \eqno (3.6)
$$

\noindent
4) for two arbitrary $t_1$-normal solutions $q_1(t)$ and $q_2(t)$ the integral
$$\ilp{t_1} Re [a(\tau) (q_1(\tau) -~ q_2(\tau))] d \tau$$ converges;

\noindent
5) for every $t_1$-normal solution $q_N(t)$ of Eq. (1.2) the equality
$$
\ilp{t_1} Re [a(\tau) (q_*(\tau) -~ q_N(\tau))] d \tau = - \infty \eqno (3.7)
$$
is valid.}

Proof. Let $q_0(t)$ be a $t_1$-regular solution of Eq. (1.2) for which $\nu_{q_0}(t_1)$ converges and $\nu_{q_0}(t) \ne 0 \ph t \ge t_1$. Then
$$
1 - \frac{1}{\nu_{q_0}(t_1)} \mu_{q_0}(t_1;t) \ne 0, \phh t \ge t_1. \eqno (3.8)
$$
Indeed otherwise if for some $t_2> t_1 \ph \nu_{q_0} = \mu_{q_0}(t_1;t_2)$ then from the equality $\nu_{q_0}(t) = \mu_{q_0}(t_1;t_2) + \nu_{q_0}(t_2)$ it follows that $\nu_{q_0}(t_2) = 0$, which contradicts our assumption. Let $q_*(t)$ be the solution of Eq. (1.2) with $q_*(t_1) = q_0(t_1) - \frac{1}{\nu_{q_0}(t_1)}$. Then by (2.7) from (3.8) it follows that $q_*(t)$ is $t_1$-regular and according to (2.11) we have
$$
\biggl|1 + \frac{1}{\nu_{q_*}(t_1)} \mu_{q_*}(t_1;t)\biggr|\biggl|1 - \frac{1}{\nu_{q_0}(t_1)} \mu_{q_0}(t_1;t)\biggr| \equiv 1, \phh t \ge t_1.
$$
From here it follows $\nu_{q_*}(t_1) = \lim\limits_{t\to +\infty}\mu_{q_*}(t_1;t) = \infty$. Then by virtue of Theorem 3.2 $q_*(t)$ is $t_1$-extremal and (3.6) is valid. Assume now Eq. (1.2) has a $t_1$-extremal solution $q_*(t)$. Show that $\nu_{q_0}(t) \ne 0, \ph t \ge t_1.$ Suppose for some $t_2 \ge t_1 \ph \nu_{q_0}(t_2) = 0$. Then obviously
$$
\lim\limits_{t \to +\infty}[1 + (q_*(t_2) - q_0(t_2))\mu_{q_0}(t_2;t)] = 1. \eqno (3.9)
$$
By (2.11) we have
$$
|1 + (q_0(t_2) - q_*(t_2))\mu_{q_*}(t_2;t)||1 + (q_*(t_2) - q_0(t_2))\mu_{q_0}(t_2;t)| \equiv 1, \phh t \ge t_2.
$$
This together wit (3.9) implies that $\mu_{q_*}(t_2;t)$ is bounded by $t$ on $[t_2,+\infty)$. Therefore  $\mu_{q_*}(t_1;t)$ is bounded by $t$ on $[t_1,+\infty)$, and according to Theorem 3.2 $q_*(t)$ is $t_1$-normal, which contradicts our assumption. The obtained contradiction shows that $\nu_{q_0}(t) \ne 0, \linebreak t \ge~ t_1$.
Let us prove (3.5). Suppose for some $t_3 \ge t_1$
$$
q_*(t_3) \ne q_0(t_3) - \frac{1}{\nu_{q_0}(t_3)}.
$$
Then there exists a finite limit
$$
\lim\limits_{t \to +\infty}[1 + (q_*(t_3) - q_0(t_3))\mu_{q_0}(t_3;t)] \ne 0. \eqno (3.10)
$$
By (2.11) we have
$$
|1 + (q_0(t_3) - q_*(t_3))\mu_{q_*}(t_3;t)| |1 + (q_*(t_3) - q_0(t_3))\mu_{q_0}(t_3;t)| \equiv 1, \ph t\ge t_3.
$$
From here and from (3.10) it follows that $\mu_{q_*}(t_3;t)$ is bounded by $t$ on $[t_3,+\infty)$. Therefore $\mu_{q_*}(t_1;t)$ is bounded by $t$ on $[t_1,+\infty)$. By virtue of Theorem 3.2 from here it follows that $q_*(t)$ is $t_1$-normal, which contradicts our assumption. The obtained contradiction proves (3.5).

Let $q(t)$ be a $t_1$-normal solution of Eq. (1.2). By (2.11) we have
$$
|1 + (q(t_1) - q_*(t_1))\mu_{q_*}(t_1;t)| |1 + (q_*(t_1) - q(t_1))\mu_q(t_1;t)| \equiv 1, \ph t\ge t_1.
$$
This together with (3.6) implies
$$
\lim\limits_{t \to +\infty} [1 + (q_*(t_1) - q(t_1))\mu_q(t_1;t)] = 0.
$$
Therefore the integrals $\nu_q(t)$ converge for all $t \ge t_1$. The inequality $\nu_q(t) \ne 0, \ph t \ge~ t_1$ follows immediately from the already proven necessary condition of existence of a \linebreak $t_1$-extremal solution of Eq. (1.2).

Let $q_1(t)$ and $q_2(t)$ be $t_1$-normal solutions of Eq. (1.2). By (2.10) we have
$$
|1 + (q_1(t_1) - q_2(t_1)) \mu_{q_2}(t_1;t)| = \exp\biggl\{\il{t_1}{t}Re [a(\tau)(q_1(\tau) - q_2(\tau))]d\tau\biggr\}, \phh t \ge t_1.
$$
From here and from the convergence of $\nu_{q_2}(t_1)$ it follows the convergence of the integral
$$
\ilp{t_1}Re [a(\tau)(q_1(\tau) - q_2(\tau))] d \tau.
$$
Let $q_N(t)$ be a $t_1$-normal solution of Eq. (1.2). By (2.10) we have
$$
|1 + (q_1(t_N) - q_*(t_1)) \mu_{q_*}(t_1;t)| = \exp\biggl\{\il{t_1}{t}Re [a(\tau)(q_*(\tau) - q_N(\tau))]d\tau\biggr\}, \phh t \ge t_1.
$$
This together with (3.6) implies (3.7). The theorem is proved.

{\bf Corollary 3.2.} {\it Let Eq. (1.2) have a $t_1$-regular solution $q_*(t)$ such that $\nu_{q_*}(t_1) = \infty.$ Then the statements 1) - 5) of Theorem 3.3 are valid.}

Proof. By Theorem 3.3 it is enough to show that Eq. (1.2) has a $t_1$-regular solution $q_0(t)$ such that $\nu_{q_0}(t_1)$ converges and $\nu_{q_0}(t) \ne 0, \ph t \ge t_1$.  Let $q_0(t)$ be a $t_1$-regular solution of Eq. (1.2), different from $q_*(t)$. In virtue of (2.11) we have
$$
|1 + (q_0(t_1) - q_*(t_1))\mu_{q_*}(t_1;t)| |1 + (q_*(t_1) - q_0(t_1))\mu_{q_0}(t_1;t)| \equiv 1, \ph t\ge t_1. \eqno (3.11)
$$
From the condition of the corollary it follows that
$$
\lim\limits_{t \to +\infty}|1 + (q_0(t_1) - q_*(t_1))\mu_{q_*}(t_1;t)| = +\infty
$$
From here and from (3.11) it follows that $q_0(t)$ is $t_1$-normal and the integral $\nu_{q_0}(t_1)$ converges. Moreover by virtue of Theorem 3.2 from the condition of the corollary it follows that $q_*(t)$ is $t_1$-extremal. Since $q_0(t)$ is an arbitrary $t_1$-regular solution of Eq. (1.2), different from $q_*(t)$ it follows that $q_*(t)$ is the unique $t_1$-extremal solution of Eq. (1.2). Then by Theorem 3.3 $\nu_{q_0}(t) \ne 0, \ph t \ge t_1.$ The corollary is proved.

Theorem 3.3 and Corollary 3.2 allow us to give the following equivalent definitions.

{\bf Definition 3.7.} {\it Eq. (1.2) is called extremal if for some $t_1 \ge t_0$ it has a $t_1$-regular solution $q(t)$ such that $\nu_q(t_1)$ converges and $\nu_q(t) \ne 0, \ph t\ge t_1$.}

{\bf Definition 3.8.} {\it Eq. (1.2) is called extremal if for some $t_1 \ge t_0$ it has a $t_1$-regular solution $q(t)$ such that $\nu_q(t_1) = \infty$.}

{\bf Example 3.3.} {\it Let $\lambda(t)$ be a quaternionic valued continuously differentiable function on $[t_0,+\infty), \ph \alpha(t) \equiv \alpha_0(t) + i \alpha_1(t), \ph \beta(t) \equiv \beta_0(t) + j \beta_1(t), \ph t\ge t_0$, where $\alpha_0(t), \ph \alpha_1(t), \ph \beta_0(t)$ and  $\beta_1(t)$ are some real-valued continuous functions on $[t_0,+\infty)$. Consider the Riccati equation
$$
q' + qa(t)q -[\lambda(t) a(t) + \alpha(t)] q - q [a(t)\lambda(t) + \beta(t)] - \lambda'(t) + \phantom{aaaaaaaaaaaaaaaaaaaaaa}
$$
$$
\phantom{aaaaaaaaaaaaaaaaaaaaa} +\lambda(t) a(t)\lambda(t) + \alpha(t)\lambda(t) + \lambda(t) \beta(t) = 0, \phh t \ge t_0. \eqno (3.12)
$$
It is not difficult to verify that $q=\lambda(t)$ is a $t_0$-regular solution of this equation and
$$
\phi_\lambda(t) = \expb{-\il{t_0}{t}\beta(\tau) d\tau}, \phh \psi_\lambda(t) = \expb{-\il{t_0}{t}\alpha(\tau)d\tau}, \ph t \ge t_0.
$$
So
$$
\nu_\lambda(t) = \ilp{t}\expb{\il{t_0}{\tau}\beta(s)}a(\tau)\expb{\il{t_0}{\tau}\alpha(s) d s}d\tau, \phh t \ge t_0.
$$
Therefore if $\nu_\lambda(t_0)$ converges and $\nu_\lambda(t) \ne 0, \ph t \ge t_t$ for some $t_1 \ge t_0$ or if $\nu_\lambda(t_0) = \infty$, then Eq. (3.12) is extremal. If  $\nu_\lambda(t_0)$ converges and $\nu_\lambda(t)$ has arbitrary large zeroes, then Eq. (1.2) is normal.}

Obviously every extremal Eq. (1.2) is sub extremal. The next example shows that not all sub extremal equations are extremal.

{\bf Example 3.4.} {\it Consider the Riccati equation
$$
q' + q (t\cos t) q = 0, \phh t \ge t_0, \ph t_0\sin t_0 + \cos t_0 = 0. \eqno (3.13)
$$
For every $\lambda \in \mathbb{H}$ the solution $q(t)$ of this equation with $q(t_0) = \lambda$ has the form
$$
q(t) = \frac{1}{1 +\lambda\il{t_0}{t}\tau\cos\tau d\tau} \lambda = \frac{1}{1 +\lambda(t\sin t + \cos t)} \lambda, \phh 1 +\lambda(t\sin t + \cos t) \ne 0.
$$
Hence every solution $q(t)$ of this equation with $q(t_0) \in \mathbb{H} \backslash  (\mathbb{R} \backslash \{0\})$ is $t_0$-regular and for $q(t_0) \in \mathbb{R}\backslash \{0\} \ph q(t)$ is not $t_0$-regular. Therefore $q_0(t) \equiv 0$ is a $t_0$-extremal solution of Eq. (3.13) and all its solutions $q(t)$ with $q(t_0) \in \mathbb{H}\backslash \mathbb{R}$ are $t_0$-normal. From here it follows that Eq. (3.13) is sub extremal. Obviously the integral
$$
\nu_{q_0}(t_0) = \ilp{t_0}t \cos t d t
$$
neither is convergent nor divergent to $\infty$. Therefore Eq. (3.13) is not extremal.}

{\bf 4. The asymptotic behavior of solutions of systems of two  first-order linear  quaternionic ordinary differential equations}. Let $a_{ml}(t), \ph m,l=1,2$ be quaternionic-valued continuous functions on $[t_0,+\infty)$. Consider the linear system
$$
\sist{\phi' = a_{11}(t)\phi + a_{12}(t)\psi,}{\psi' = a_{21}(t)\phi + a_{22}(t)\psi, \ph t \ge t_0} \eqno (4.1)
$$
and the quaternionic Riccati equation
$$
q' + q a_{12}(t) q + q a_{11}(t) - a_{22}(t) q - a_{21}(t) = 0, \phh t \ge t_). \eqno (4.2)
$$
It is not difficult to verify that the solutions $q(t)$ of Eq. (4.2), existing on some interval $[t_1,t_2) \ph (t_0 \le t_1 < t_2 \le +\infty)$ are connected with solutions $(\phi(t), \psi(t))$ of the system (4.1) by relations
$$
\phi'(t) = [a_{12}(t) q(t) + a_{11}(t)] \phi(t), \phh \psi(t) = q(t)\phi(t), \phh t \in [t_1,t_2). \eqno (4.3)
$$
From here it follows
$$
\widehat{\phi(t)}' = [\widehat{a_{12}(t)} \widehat{q(t)} + \widehat{a_{11}(t)}]\widehat{\phi(t)}, \ph t\in [t_1,t_2).
$$
By Liouville's formula from here we obtain
$$
det\hskip 3pt \widehat{\phi(t)} = det \hskip 3pt \widehat{\phi(t_1)} \exp\biggl\{\il{t_1}{t} tr \hskip 3pt [\widehat{a_{12}(\tau)}\widehat{q(t)} + \widehat{a_{11}(\tau)}]d \tau\biggr\}, \phh t \in[t_1;t_2).
]$$
By virtue of Lemma 2.1 from here it follows
$$
|\phi(t)| = |\phi(t_1)|\expb{\il{t_1}{t} Re \hskip 3pt [a_{12}(\tau) q(\tau) + a_{11}(\tau)] d \tau}, \ph t \in [t_1,t_2). \eqno (4.4)
$$
So if $\phi(t_1) \ne 0$, then
$$
\phi(t) \ne 0, \phh t \in [t_1,t_2). \eqno (4.5)
$$

{\bf Remark 4.1.} {\it It can be shown that if for a solution $(\phi(t), \psi(t))$ of the system (4.1) the function $\phi(t)$ does not vanish on $[t_1,t_2)$ then $q(t) = \psi(t)\phi^{-1}(t), \ph t\in [t_1,t_2)$ is a solution of Eq. (4.2) on $[t_1,t_2)$.}

{\bf Definition 4.1.} {\it A solution $(\phi(t), \psi(t))$ of the system (4.1) is called $t_1$-regular ($t_1 \ge t_0$)
if $\phi(t) \ne 0, \ph t \ge t_1$.}

{\bf Definition 4.2.} {\it A $t_1$-regular ($t_1\ge t_0$) solution $(\phi(t), \psi(t))$ of the system (4.1) is called principal (non principal) if $q(t)\equiv \psi(t)\phi^{-1}(t), \ph t \ge t_1$ is a $t_1$-extremal ($t_1$-normal)
solution of Eq. (4.2).}

{\bf Definition 4.3.} {\it The system (4.1) is called regular if it has at least one $t_1$-regular solution for some $t_1 \ge t_0$.}

{\bf Remark 4.2.} {\it It follows from (4.5) and Remark 4.1 that the system (4.1) has a $t_1$-regular solution for some $t_1 \ge t_0$ if and only if Eq. (4.2) has a $t_1$-regular solution.}

{\bf Remark 4.3.} {\it If $(\phi(t), \psi(t))$ is a solution of the system (4.1) then for every $\lambda \in~ \mathbb{H} \linebreak (\phi(t)\lambda, \psi(t)\lambda)$ is also a solution of the system (4.1), but $(\lambda\phi(t), \lambda \psi(t))$ may not be a solution of the system (4.1). For example $(e^{it}, e^{kt}), \ph t \ge t_0$ is a solution of the system
$$
\sist{\phi' = i \phi,}{\psi' = k\psi, \ph t \ge t_0}
$$
but $(j e^{it}, j e^{kt}), \ph t \ge t_0$ is not a solution of this system.}

{\bf Definition 4.4.} {\it The solutions $(\phi_m(t), \psi_m(t)), \ph m=1,2$ are called linearly dependent if there exists $\lambda \in \mathbb{H} \backslash \{0\}$ such that $\phi_2(t) = \phi_1(t)\lambda, \ph \psi_2(t) = \psi_1(t)\lambda$, otherwise they are called linearly independent.}

{\bf Remark 4.4.} {\it It follows from Theorem 3.1 and Remark 4.1 that if the system (4.1) has a $t_1$-regular solution $(\phi(t), \psi(t))$, then it has also another $t_1$-regular solution, linearly independent of $(\phi(t), \psi(t)).$}

{\bf Definition 4.5.} {\it The regular system  (4.1)  is called normal (irreconcilable, sub extremal, super extremal, extremal)  if Eq. (4.2) is normal (irreconcilable, sub extremal, super\linebreak extremal, extremal).}

Hereafter every $t_1$-regular solution of the system (4.1) we will just call a regular solution of the system (4.1). On the basis of  (4.4) from Corollary 3.1 we immediately get.

{\bf Theorem 4.1.} {\it The following statements are valid:

\noindent
I) if the system (4.1) is normal then for its two regular solutions $(\phi_m(t), \psi_m(t)), \ph m=1,2$ the inequalities
$$
\limsup\limits_{t \to +\infty} \frac{|\phi_1(t)|}{|\phi_2(t)|} < +\infty, \phh \limsup\limits_{t \to +\infty} \frac{|\phi_2(t)|}{|\phi_1(t)|} < +\infty
$$
are valid;

\noindent
II)  if the system (4.1) is irreconcilable   then for its two arbitrary linearly independent regular solutions $(\phi_m(t), \psi_m(t)), \ph m=1,2$ the equalities
$$
\limsup\limits_{t \to +\infty} \frac{|\phi_1(t)|}{|\phi_2(t)|} = \limsup\limits_{t \to +\infty} \frac{|\phi_2(t)|}{|\phi_1(t)|} = +\infty
$$
are valid;

\noindent
III) If the system (4.1) is sub extremal then there exists a regular solution $(\phi_*(t), \psi_*(t))$ of (4.1) such that for every regular solutions $(\phi_m(t), \psi_m(t)), \ph m=1,2$ of (4.1) linearly independent of $(\phi_*(t), \psi_*(t))$ the relations
$$
\limsup\limits_{t \to +\infty} \frac{|\phi_*(t)|}{|\phi_1(t)|} < +\infty, \phh \liminf\limits_{t \to +\infty} \frac{|\phi_*(t)|}{|\phi_1(t)|} = 0,
$$
$$
\limsup\limits_{t \to +\infty} \frac{|\phi_1(t)|}{|\phi_2(t)|} < +\infty, \phh \limsup\limits_{t \to +\infty} \frac{|\phi_2(t)|}{|\phi_1(t)|} < +\infty
$$
are valid;

\noindent
IV) if the system (4.1) is super extremal then there exist two regular solutions $(\phi_*(t), \psi_*(t))$ and  $(\phi^*(t), \psi^*(t))$ of (4.1) such that
$$
\limsup\limits_{t \to +\infty} \frac{|\phi_*(t)|}{|\phi^*(t)|} = \limsup\limits_{t \to +\infty} \frac{|\phi^*(t)|}{|\phi_*(t)|} = +\infty
$$
and for all  two  arbitrary solutions $(\phi_m(t), \psi_m(t)), \ph m=1,2$ of (4.1) linearly independent of each   $(\phi_*(t), \psi_*(t))$ and  $(\phi^*(t), \psi^*(t))$  the following relations are valid
$$
\limsup\limits_{t \to +\infty} \frac{|\phi_1(t)|}{|\phi_2(t)|} < +\infty, \phh \limsup\limits_{t \to +\infty} \frac{|\phi_2(t)|}{|\phi_1(t)|} < +\infty,
$$
$$
\limsup\limits_{t \to +\infty} \frac{|\phi_*(t)|}{|\phi_m(t)|} < +\infty, \phh \limsup\limits_{t \to +\infty} \frac{|\phi^*(t)|}{|\phi_m(t)|} < +\infty,
$$
$$
 \liminf\limits_{t \to +\infty} \frac{|\phi_*(t)|}{|\phi_m(t)|}=
\liminf\limits_{t \to +\infty} \frac{|\phi^*(t)|}{|\phi_m(t)|}=0, \ph m=1,2.
$$
}~
\phantom{aaaaaaaaaaaaaaaaaaaaaaaaaaaaaaaaaaaaaaaaaaaaaaaaaaaaaaaaaaaaaaaaaaaa}$\blacksquare$

Theorem 4.1 shows that in the normal case of the system (4.1) all regular solutions of (4.1) are asymptotically equivalent.  This case differs from the other cases by the scarcity of asymptotic behavior patterns at $+\infty$ of the solutions of the system (4.1). In the supercritical case of (4.1) we have "the richest" (among the other cases) variety of asymptotic behavior pattern at $+\infty$ of regular solutions of the system (4.1)

Let
$$
a_{12}(t) = a_0(t) + i a_1(t) + j a_2(t) + k a_3(t), \ph -a_{22}(t) = b_0(t) + i b_1(t) + j b_2(t) + k b_3(t),
$$
$$
a_{11}(t) = c_0(t) + i c_1(t) + j c_2(t) + k c_3(t), \ph -a_{21}(t) = d_0(t) + i d_1(t) + j d_2(t) + k d_3(t).
$$
where $a_m(t), \ph b_m(t), \ph c_m(t)$ and $d_m(t), \ph m=\overline{0,3}$ are real-valued continuous functions on $[t_0,+\infty)$. Set:
$$
p_{0,m}(t)\equiv b_m(t) + c_m(t), \ph m=\overline{1,3}
$$
$$
p_{11}(t) \equiv b_1(t) + c_1(t), \phh p_{12}(t) \equiv b_2(t) - c_2(t),
$$
$$
p_{13}(t) \equiv b_3(t) - c_3(t), \phh p_{21}(t) \equiv b_1(t) - c_1(t),
$$
$$
p_{22}(t) \equiv b_2(t) + c_2(t), \phh p_{23}(t) \equiv b_3(t) - c_3(t),
$$
$$
p_{3m}(t)\equiv b_m(t) - c_m(t), \ph m=\overline{1,3}, \ph t \ge t_0,
$$
$$
D_0(t) \equiv \sist{\sum\limits_{m=1}^{3}p_{0m}^2(t) + 4 a_0(t) d_0(t), \ph \mbox{if} \ph a_0(t) \ne 0,}{4d_0(t) \phh \mbox{if} \ph a_0(t) =0,} \phantom{aaaaaaaaaaaaaaaaaaaaaaaaaaaaaaaaaaaa}
$$
$$
\phantom{aaaaaaaaa}D_n(t) \equiv \sist{\sum\limits_{m=1}^{3}p_{nm}^2(t) - 4 a_n(t) d_n(t), \ph \mbox{if} \ph a_n(t) \ne 0,}{-4d_n(t) \phh \mbox{if} \ph a_n(t) =0,} \phh n-\overline{1,3}, \phh t \ge t_0.
$$

Let $\mathfrak{S}$ be a non empty subset of the set $\{0, 1,  2,  3\}$ and let $\mathfrak{D}$ be its complement i. e. $\mathfrak{D} = \{0, 1,  2,  3\} \backslash \mathfrak{S}.$

{\bf Theorem 4.2}. {\it Let the conditions

\noindent
$\alpha)  \ph  a_n(t) \ge 0, \ph t \ge t_0, \ph n\in \mathfrak{S}$ and if $a_n(t) = 0$ then $p_{nm}(t) = 0, \ph m\in \mathfrak{S}, \ph a_n(t) \equiv~ 0, \linebreak n\in~ \mathfrak{D}, \ph D_n(t) \le~ 0, \ph t \ge t_0, \ph n = \overline{0,3}$;

\noindent
$\beta) \ph \ilp{t_0}|a_{12}(\tau)|\expb{\il{t_0}{t} \Bigl[\re a_{22}(s) - \re a_{11}(s)\Bigr] d s} d \tau < +\infty$.

\noindent
be satisfied. Then the following statements are valid:

\noindent
1) \ph the system (4.1) is or else normal or else extremal:

\noindent
2) \ph for all $T$-regular ($T \ge t_0$) non principal solutions $(\phi(t), \psi(t))$ of the system (4.1) the integral
$$
\ilp{T}\frac{|a_{12}(\tau)|}{|\phi(\tau)|^2}\expb{\il{T}{\tau}\Bigl[\re a_{11}(s) + \re a_{22}(s)\Bigr]d s} d \tau
$$
converges;

\noindent
3) \ph  if the system (4.1) is extremal, then:

\noindent
$3_1$) \ph for its unique (up to arbitrary right multiplier) principal solution $(\phi_*(t), \psi_*(t))$ the equality
$$
\ilp{T_*}\frac{|a_{12}(\tau)|}{|\phi_*(\tau)|^2}\expb{\il{T_*}{\tau}[\re a_{11}(s) + \re  a_{22}(s)]d s}d\tau = +\infty; \eqno (4.6)
$$
is valid, where $T_* \ge t_0$ such that $\phi_*(t)\ne 0, \ph t \ge T_*$;

\noindent
$3_2$) \ph for all non principal solutions $(\phi(t), \psi(t))$ of the system (4.1) the equality
$$
\liml\frac{|\phi_*(t)|}{|\phi(t)||} = 0  \eqno (4.7)
$$
is valid;

\noindent
$3_3$) \ph for two arbitrary non principal solutions $(\phi_m(t), \psi_m(t)), \ph m=1,2$ of the system (4.1) the relation
$$
\liml\frac{|\phi_1(t)|}{|\phi_2(t)||} = c \ne 0  \eqno (4.8)
$$
is valid.}

To prove this theorem we need in the following result from [7] (see [7, Theorem 3.1])

{\bf Theorem 4.3}. {\it Let the conditions $\alpha$) of Theorem 4.2 be satisfied. Then for all $\gamma_n \ge~ 0, \linebreak n\in \mathfrak{S}, \ph \gamma_n \in (-\infty, +\infty), \ph n \in \mathfrak{D}$ Eq. (4.2) has a solution $\mathfrak{q}_0(t) = \mathfrak{q}_{0, 0}(t) - i \mathfrak{q}_{0, 1}(t) - j \mathfrak{q}_{0, 2}(t) - k \mathfrak{q}_{0, 3}(t)$ on $[t_0, +\infty)$ with $\mathfrak{q}_{0, n}(t_0) = \gamma_n, \ph n= \overline{0, 3}$ and  $\mathfrak{q}_{0, n}(t) \ge 0, \ph  n \in \mathfrak{S}, \ph t \ge t_0$.}

\phantom{aaaaaaaaaaaaaaaaaaaaaaaaaaaaaaaaaaaaaaaaaaaaaaaaaa} $\blacksquare$

{\bf Proof of Theorem 4.2.} Let $q_0(t)$ be the solution of Eq. (4.2) with $q_0(t) = 0$. In virtue of Theorem 4.3 it follows from the conditions $\alpha$) of the theorem that $q_0(t)$ is $t_0$-regular and
$$
\re [a_{12}(t) q_0(t)] \ge 0, \phh t \ge t_0. \eqno (4.9)
$$
Consider the integral
$$
\widetilde{\nu}_{q_0}(t) \equiv \ilp{t}\phi^{-1}_{q_0}(\tau) a_{12}(\tau)\psi^{-1}_{q_0}(\tau) d\tau, \phh t \ge t_0,
$$
where $\phi_{q_0}(t)$ and $\psi_{q_0}(t)$ are the solutions of the linear equations
$$
\phi' = [a_{12}(t) q_0(t) + a_{11}(t)]\phi, \phh t \ge t_0, \eqno (4.10)
$$
$$
\psi' = \psi[q_0(t)a_{12}(t) - a_{22}(t)], \phh t \ge t_0
$$
respectively with $\phi_{q_0}(t_0) = \psi_{q_0}(t_0) =1.$ By  (2.7) and (2.8) we have respectively
$$
|\phi_{q_0}(t)| = \expb{\il{t_0}{t}\re [a_{12}(\tau)q_0(\tau) + a_{11}(\tau)]}, \eqno (4.11)
$$
$$
|\psi_{q_0}(t)| = \expb{\il{t_0}{t}\re [a_{12}(\tau)q_0(\tau) - a_{22}(\tau)]}, \ph t \ge t_0.
$$
Hence,
$$
|\widetilde{\nu}_{q_0}(t)| \le \ilp{t}\frac{| a_{12}(\tau)|}{|\phi_{q_0}(\tau)||\psi_{q_0}(\tau)|} d\tau= \phantom{aaaaaaaaaaaaaaaaaaaaaaaaaaaaaaaaaaaaaaaaaaaaaaaaaaaa}
$$
$$
\phantom{aaaaaaaaa}=\ilp{t}|a_{12}(\tau)|\expb{-\il{t_0}{\tau}[2\re a_{12}(s) q_0(s) + \re a_{11}(s) - \re a_{22}(s)]d s}d\tau, \ph t \ge t_0.
$$
This together with (4.9) and $\beta$) implies that
$$
|\widetilde{\nu}_{q_0}(t)| \le \ilp{t}|a_{12}(\tau)|\expb{\il{t_0}{\tau}[\re a_{22}(s) - \re a_{11}(s)]d s}d \tau < +\infty \phh t \ge t_0. \eqno (4.12)
$$
It follows from here that the integrals $\widetilde{\nu}_{q_0}(t), \ph t \ge t_0$ converge.
Two cases are possible:

\noindent
$a) \ph \widetilde{\nu}_{q_0}(t)$ has arbitrary large zeroes;

\noindent
$b) \ph \widetilde{\nu}_{q_0}(t) \ne 0, \ph t \ge T_0$ for some $T_0 \ge t_0$.

\noindent
Then by Theorem 3.3 the system (4.1) is or else normal (in the case $a$)) or else extremal (in the case $b$)). The statement 1) of the theorem is proved.
Let $(\phi_0(t),\psi_0(t))$ be the solution of the system (4.1) with $\phi_0(t_0) = 1, \ph \psi_0(t_0) =~0$. Then by (4.3) $\phi_0(t)$ is a solution of Eq.~ (4.10). So $\phi_0(t)$ coincides with $\phi_{q_0}(t)$. Therefore  from $\beta$),  (4.9) and (4.11) it follows
$$
\ilp{t}\frac{|a_{12}(\tau)|}{|\phi_0(\tau)|^2}\expb{\il{t_0}{\tau}[\re a_{11}(s) + \re a_{22}(s)]d s} d\tau \le \phantom{aaaaaaaaaaaaaaaaaaaaaaaaaaaaaaa}
$$
$$
\phantom{aaaaaaaaaaa}\le \ilp{t}|a_{12}(\tau)|\expb{\il{t_0}{\tau}[\re a_{22}(s) - \re a_{11}(s)]d s}d\tau < +\infty, \ph t \ge t_0. \eqno (4.13)
$$
Let $(\phi(t), \psi(t))$ be a $T$-regular ($T\ge t_0$) non principal solution of the system (4.1). Then $q(t)\equiv \psi(t)\phi^{-1}(t), \ph t \ge T$ is a $T$-normal solution of Eq. (4.2). It follows from (4.12) that $\mu_{q_0}(T;t)$ is bounded on $[T,+\infty)$. Hence, according to the statement 1) of Corollary 3.1 we have
$$
\Supl{t\ge T} \biggl|\il{T}{t}\re [a_{12}(\tau)(q_0(\tau) - q(\tau))]d\tau\biggr| < +\infty.
$$
This together with (4.10) implies
$$
\ilp{T}\frac{|a_{12}(\tau)|}{|\phi(\tau)|^2}\expb{\il{T}{\tau}[\re a_{11}(s) + \re a_{22}(s)]d s} d\tau =
$$
$$
=\ilp{T}\frac{|a_{12}(\tau)|}{|\phi_0(\tau)|^2}\expb{\il{T}{\tau}[\re a_{11}(s) + \re a_{22}(s)]d s}\expb{2\il{T}{\tau}\re [a_{12}(s)(q_0(s) - q(s))]d s} d\tau\le
$$
$$
\phantom{aaaaaaaaaaaaaaaaaaaaaa}\le M \ilp{T}\frac{|a_{12}(\tau)|}{|\phi_0(\tau)|^2}\expb{\il{t_0}{\tau}[\re a_{11}(s) + \re a_{22}(s)]d s}d\tau < +\infty,
$$
where
$$
M\equiv  \expb{-\il{t_0}{T}[\re a_{11}(s) + \re a_{22}(s)]d s} \expb{2\Supl{t\ge T}\Bigl|\il{t_0}{\tau}\re [a_{12}(s)(q_(s) - q(s))]d s\Bigr|}<~ +\infty.
$$
The statement 2) of  the theorem is proved.  Assume the system (4.1) is extremal. Then Eq. (4.2) has the unique extremal solution $q_*(t)$. Let $q_*(t)$ be $T_*$-regular for some $T_*\ge t_0$ and let $(\phi_*(t), \psi_*(t))$ be the solution of the system (4.1) with $\phi_*(T_*) =1, \ph \psi_*(T_*) =q_*(T_*)$. Then by (4.3) $(\phi_*(t), \psi_*(t))$ is the unique (up to arbitrary right multiplier) principal solution of the system (4.1) and $\phi_*(t)$ is a solution of the linear equation
$$
\phi' = [a_{12}(t)q_*(t) + a_{11}(t)]\phi, \phh t \ge T_*. \eqno (4.14)
$$
 Consider the integral
$$
\widetilde{\nu}_{q_*}(T_*) \equiv\ilp{T_*}\phi_{q_*}^{-1}(\tau)a_{12}(\tau)\psi_{q_*}^{-1}(\tau) d \tau,
$$
where $\phi_{q_*}(t)$ and $\psi_{q_*}(t)$ are the solutions of Eq. (4.14) and  the  equation
$$
\psi' = \psi[q_*(t)a_{12}(t) - a_{22}(t)], \phh t \ge T_*
$$
respectively with $\phi_{q_*}(T_*) = \psi_{q_*}(T_*) = 1$. Since $q_*(t)$ is extremal in virtue of Theorem 3.3 we have
$$
\widetilde{\nu}_{q_*}(T_*) = \infty. \eqno (4.15)
$$
By (2.7) and (2.8) we have respectively
$$
|\phi_{q_*}(t)| = \expb{\il{T_*}{t}\re [a_{12}(\tau) q_*(\tau) + a_{11}(\tau)]d\tau}, \phh t \ge T_*,
$$
$$
|\psi_{q_*}(t)| = \expb{\il{T_*}{t}\re [a_{12}(\tau) q_*(\tau) - a_{22}(\tau)]d\tau}, \phh t \ge T_*.
$$
Therefore
$$
|\psi_{q_*}(t)| = |\phi_{q_*}(t)|\expb{-\il{T_*}{t}\re [a_{11}(\tau) + a_{22}(\tau)]d\tau}, \phh t \ge T_*. \eqno (4.16)
$$
Obviously $\phi_*(t) = \phi_{q_*}(t), \ph t \ge T_*$. This together with  (4.16) implies
$$
|\widetilde{\nu}_{q_*}(T_*)| \le \ilp{T_*}\frac{|a_{12}(\tau)|}{|\phi_*(\tau)|^2}\expb{\il{T_*}{\tau}\re[a_{11}(s) + a_{22}(s)]d s}d\tau.
$$
From here and from (4.15) it follows (4.6). Let $(\phi(t), \psi(t))$ be a non principal solution of the system (4.1). Without loss of generality we may take that $(\phi(t), \psi(t))$ is $T_*$-regular. Then $q(t)\equiv \psi(t) \phi^{-1}(t), \ph t \ge T_*$ is a $T_*$-normal solution of Eq. (4.2). By (3.7) from here it follows
$$
\ilp{T_*}\re [a_{12}(\tau) (q_*(\tau) - q(\tau))] d\tau = -\infty.
$$
By (2.7) from here we obtain (4.7):
$$
\liml\frac{|\phi_*(t)|}{|\phi(t)|} = \liml\expb{\il{T_*}{t}\re [a_{12}(\tau) (q_*(\tau) - q(\tau))] d\tau} = 0.
$$
Let $(\phi_m(t), \psi_m(t)), \ph m=1,2$  be non principal $T$-regular ($T \ge t_0$) solutions of the system (4.1). By (4.3) $q_m(t) = \psi_m(t)\phi_m^{-1}(t), \ph t \ge~ T, \ph m=1,2$ are $T$-normal solutions of Eq. (4.2). Then according to the statement 4) of Theorem 3.3 the integral
$$
\ilp{T}\re [a_{12}(\tau) (q_1(\tau) - q_2(\tau))] d\tau
$$
converges.
 By (2.7) from here it follows (4.8). The theorem is proved.

{\bf Remark 4.1.} {\it From the estimate (4.12) is seen that if $supp \hskip 1.5pt a_{12}(t)$ is bounded, then $\widetilde{\nu}_{q_0}(t)$ has arbitrary large zeroes. Hence in this case under the conditions of Theorem 4.2 the system is normal. If $supp \hskip 1.5pt a_{12}(t)$ is unbounded and the coefficients of the system (4.1) are real-valued,  then it is not difficult to verify that under the conditions of Theorem 4.1  $\widetilde{\nu}_{q_0}(t)\ne 0, \ph t \ge t_0$. So in this case (4.1) is extremal.}

\vskip 20pt

\centerline{\bf References}

\vskip 20pt

\noindent
1. P. Wilzinski,  Quaternionic-valued differential equations. The Riccati equations. \linebreak\phantom{aa}  Journal of Differential Equations, vol. 247. pp. 2167 - 2187, 2009.

\noindent
2. J. D. Gibbon, D. D. Holm, R. M. Kerr and I. Roulstone. Quaternions and periodic \linebreak\phantom{aa}  dynamics in the Euler fluid equations. Nonlinearity, vol. 19, pp. 1962 - 1983, 2006.

\noindent
3. H. Zoladek, Classification of diffeomorphisms of $\mathbb{S}^4$ induced by quaternionic Riccati \linebreak\phantom{aa}  equations with periodic coefficients. Topological methods in Nonlinear Analysis. Journal \linebreak\phantom{aa}  of the Juliusz Shauder Center, vol. 33. pp. 205 - 2015, 2009.

\noindent
4. V. Cristioano and F. Smaraadase. An Exact Mapping from Navier-Stocks Equation to \linebreak\phantom{aa} Schrodinger Equation via Riccati equation. Progress in Physics, vol. 1. pp. 38, 39, 2008.

\noindent
5. K. Leschke and K. Morya. Application of Quaternionic Holomorphic Geometry to \linebreak\phantom{aa} minimal surfaces. Complex manifolds, vol. 3, pp. 282 - 300, 2006.

\noindent
6. J. Campos, J. Mavhin. Periodic solutions of quaternionic-valued ordinary differential \linebreak\phantom{aa} equations. Annali di Mathematica, vol. 185, pp. 109 - 127, 2006.

\noindent
7. G. A. Grigorian, Global solvability criteria for quaternionic Riccati equations. Archivum \linebreak\phantom{aa}   Mathematicum. In print.

\noindent
8. G. A. Grigorian, On some properties of solutions of the Riccati equation.
Izvestiya  NAS \linebreak\phantom{aa} of Armenia, vol. 42, $N^\circ$ 4, 2007, pp. 11 - 26.

\noindent
9. G. A. Grigorian, On the Stability of Systems of Two First-Order Linear Ordinary  \linebreak\phantom{aa}
Differential Equations, Differ. Uravn., 2015, vol. 51, no. 3, pp. 283 - 292.

\noindent
10. G. A. Grigorian. Necessary Conditions and a Test for the Stability of a System of Two  \linebreak\phantom{aa}
Linear Ordinary Differential Equations of the First Order. Difer. Uravn., 2016,
Vol.~ 52, \linebreak\phantom{aa}  No. 3, pp. 292 - 300.

\noindent
11.  G. A. Grigorian, On one oscillatory criterion for the second order linear
 ordinary \linebreak \phantom{a} differential equations. Opuscula Math. 36, no. 5 (2016), 589–601. \linebreak \phantom{aa}
   http://dx.doi.org/10.7494/OpMath.2016.36.5.589

\noindent
12. G. A. Grigorian. Some properties of the solutions of third order linear ordinary \linebreak \phantom{a}    differential  equations.  Rocky Mountain Journal of Mathematics, vol. 46, no. 1,  2016, \linebreak \phantom{a}  pp. 147 - 161.

\noindent
13. G. A. Grigorian.    Oscillatory criteria for the second order linear ordinary differential \linebreak \phantom{a} equations.
Math. Slovaca 69 (2019), No. xx,  1- 14.

\noindent
14. A. I. Egorov. Riccati equations. Moskow, Fizmatlit,  2001.

\end{document}